\documentclass[12pt]{amsart}

\newtheorem{theorem}{Theorem}
\theoremstyle{remark}
\newtheorem*{question}{Open question}
\newtheorem{remark}{Remark}


\begin{document}

\title{Bohr's power series theorem in several variables}

\author{Harold P. Boas}

\address{Department of Mathematics, Texas A\&M
University, 
College Station,
TX 77843--3368}

\thanks{The first author's research was supported in part by NSF
grant number DMS 9500916 and in part at 
the Mathematical Sciences Research Institute 
by NSF grant number DMS 9022140.}

\email{boas@math.tamu.edu}

\author{Dmitry Khavinson}
\address{Department of Mathematical Sciences, University of
Arkansas, Fayetteville, AR 72701}
\email{dmitry@comp.uark.edu}


\begin{abstract}
Generalizing a classical one-variable theorem of Bohr, we show
that if an $n$-variable power series has modulus less than~$1$
in the unit polydisc, then the sum of the moduli of the terms is
less than~$1$ in the polydisc of radius $1/(3\sqrt n\,)$.
\end{abstract}

\maketitle

How large can the sum of the moduli of the terms of a convergent
power series be? Harald Bohr addressed this question in~1914 with
the following
remarkable result on power series in one complex variable.

\begin{theorem}[Bohr]
Suppose that the power series $\sum_{k=0}^\infty c_k z^k$
converges for $z$ in the  unit disk, and $|\sum_{k=0}^\infty c_k
z^k|< 1$ when $|z|< 1$.  Then $\sum_{k=0}^\infty |c_k z^k|< 1$
when $|z|< 1/3$. Moreover, the radius $1/3$ is the best possible.
\end{theorem}

Bohr's paper \cite{Bohr}, compiled by G.~H.~Hardy from
correspondence, indicates that Bohr initially obtained the radius
$1/6$, but this was quickly improved to the sharp result by
M.~Riesz, I.~Schur, and N.~Wiener, independently.  Bohr's paper
presents both his own proof and Wiener's. Some years later,
S.~Sidon gave a different proof \cite{Sidon}, which was
subsequently rediscovered by M.~Tomi\'c \cite{Tomic}.

In this note, we formulate a version of Bohr's theorem in higher
dimensions. We write an $n$-variable power series $\sum_\alpha
c_\alpha z^\alpha$ using the standard multi-index notation:
$\alpha$~denotes an $n$-tuple
$(\alpha_1,\alpha_2,\dots,\alpha_n)$ of nonnegative integers,
$|\alpha|$~denotes the sum $\alpha_1+\dots+\alpha_n$ of its
components, 
$\alpha!$~denotes the product $\alpha_1!\alpha_2!\dots\alpha_n!$
of the factorials of its components,
$z$~denotes an $n$-tuple $(z_1,\dots,z_n)$ of complex
numbers, and $z^\alpha$~denotes the product
$z_1^{\alpha_1}z_2^{\alpha_2}\dots z_n^{\alpha_n}$.

Let $K_n$ denote the $n$-dimensional \emph{Bohr radius}: the
largest number such that if $\sum_\alpha c_\alpha z^\alpha$ converges in
the unit polydisc $\{(z_1,\dots,z_n):\max_{1\le j\le
n}|z_j|<1\}$, and if $|\sum_\alpha c_\alpha z^\alpha|<1$ in the unit
polydisc, then $\sum_\alpha|c_\alpha z^\alpha|<1$ when
$\max_{1\le j\le n}|z_j|<K_n$. 

It is evident from Bohr's one-dimensional result that $K_n\le
1/3$ for every~$n$. Surprisingly, $K_n$~tends to~$0$ when
$n$~increases. Our result implies that the decay rate of~$K_n$ is
essentially $1/\sqrt n$.

\begin{theorem}
\label{bohr}
When $n>1$, the  $n$-dimensional
Bohr radius~$K_n$ satisfies
\begin{equation*}
\frac{(1/3)}{\sqrt n} < K_n < \frac{2\sqrt{\log n}}{\sqrt n}.
\end{equation*}
\end{theorem}

Although the theorem gives tight control on the Bohr radius,
there is a small amount of slack that we have not been able to
remove.

\begin{question}
What is the exact value of the Bohr radius~$K_n$ when $n>1$?
\end{question}

There is an analogous concept of a Bohr radius for the unit ball,
or more generally for domains of convergence for power series,
that is, for logarithmically convex Reinhardt domains (for
definitions, see, for example, section~2.3 of~\cite{Krantz}). It
is easy to see that the polydisc constant~$K_n$ is a universal
lower bound for the Bohr radius of every multi-circular domain
(whether logarithmically convex or not).

\begin{theorem}
\label{multi}
Suppose that the $n$-variable power series $\sum_{\alpha}
c_\alpha z^\alpha$ converges in a complete Reinhardt domain~$G$,
and $|\sum_{\alpha} c_\alpha z^\alpha|<1$ when $z\in G$. Then
$\sum_{\alpha}|c_\alpha z^\alpha|< 1$ when $z$~is in the scaled
domain $K_n\cdot G$, where $K_n$~is the Bohr radius for the unit polydisc.
\end{theorem}

Indeed, a linear change of variables shows that
Theorem~\ref{multi} holds when $G$~is any bounded polydisc (with
its $n$~radii not necessarily equal to each other). By
definition, a complete Reinhardt domain~$G$ is a union of
polydiscs centered at the origin, so $G$~inherits the conclusion
of the theorem from the polydisc case.

\begin{proof}[Proof of Theorem~\ref{bohr}]
Concerning the lower bound, we will prove more than is stated in
the theorem: namely, if $f(z)=\sum_{\alpha}c_\alpha z^\alpha$ is
an analytic function of modulus less than~$1$ in the unit
polydisc, then $\sum_\alpha |c_\alpha z^\alpha|<1$ in the ball of
radius~$1/3$ centered at the origin. This ball evidently contains
the polydisc $\{z:\max_{1\le j\le n}|z_j|<1/(3\sqrt n\,)\}$,
whence $K_n\ge 1/(3\sqrt n\,)$. (It will become apparent below
why strict inequality obtains.)

We begin by establishing bounds on the coefficients~$c_\alpha$.
By Cauchy's estimate, $|c_\alpha|\le1$ for every~$\alpha$. We
improve this estimate for $\alpha$ different from the
$0$~multi-index by an argument analogous to one used by Wiener in
the single-variable case.

Let $\omega$~denote a primitive $k$th root of unity, and let $g$
be defined by $g(z)=k^{-1}\sum_{j=1}^k f(\omega^j z)$. The
modulus of~$g$ is again less than~$1$ in the unit polydisc, and
the Taylor series of~$g$ begins $c_0+\sum_{|\alpha|=k}c_\alpha
z^\alpha+\dots$ (because the nonconstant terms of homogeneity
less than~$k$ average out).  Next define~$h$ via $h(z)=(
g(z)-c_0)/(1-\bar c_0 g(z))$.  The modulus of~$h$ is less
than~$1$ in the unit polydisc, and the Taylor series of~$h$
begins $\sum_{|\alpha|=k}b_\alpha z^\alpha+\dots$, where the
coefficients of the leading terms (the ones of homogeneity~$k$)
satisfy $b_\alpha=c_\alpha/(1-|c_0|^2)$. Since the modulus of~$h$
does not exceed~$1$ in the unit polydisc, neither does its
$L^2$~norm on the unit torus.  (We normalize Lebesgue measure on
the unit torus to have total mass~$1$.) The monomials~$z^\alpha$ are
orthonormal on the unit torus, so we have (in particular) that
$\sum_{|\alpha|=k}|b_\alpha|^2\le 1$, whence
$(\sum_{|\alpha|=k}|c_\alpha|^2)^{1/2} \le (1-|c_0|^2)$.

The Cauchy-Schwarz inequality now implies that
\begin{equation}
\label{eqn:est}
\begin{aligned}
\sum_{\alpha}|c_\alpha z^\alpha|&= |c_0|+\sum_{k=1}^\infty
\sum_{|\alpha|=k} |c_\alpha z^\alpha|\\
& \le |c_0|+(1-|c_0|^2)\sum_{k=1}^\infty \biggl(\sum_{|\alpha|=k}
|z^\alpha|^2 \biggr)^{1/2}.
\end{aligned}
\end{equation}
But $\sum_{|\alpha|=k}|z^\alpha|^2 \le
(\sum_{j=1}^n|z_j|^2)^k$, so if $z$~lies in the ball
of radius~$1/3$, then 
\begin{equation}
\label{eqn:main}
\sum_{\alpha}|c_\alpha z^\alpha| \le |c_0| + (1-|c_0|^2)
\sum_{k=1}^\infty \frac{1}{3^k} = |c_0|+\tfrac12 (1-|c_0|^2).
\end{equation}
The right-hand side of~(\ref{eqn:main}) does not exceed~$1$,
whatever the value of~$c_0$.
Thus $\sum_\alpha|c_\alpha z^\alpha|<1$ in the ball of
radius~$1/3$, and so $K_n\ge 1/(3\sqrt n\,)$. 

If $k>1$ and at least two of the coordinates~$z_j$ are nonzero,
then $\sum_{|\alpha|=k} |z^\alpha|^2$ is strictly less than
$(\sum_{j=1}^n|z_j|^2)^k$. This means that the set of~$z$
for which (\ref{eqn:main})~holds is a logarithmically convex
Reinhardt domain slightly fatter than the ball of
radius~$1/3$. Hence $K_n$~strictly exceeds $1/(3\sqrt n\,)$ when
$n>1$.

We now turn to the right-hand inequality in the
theorem. According to the theory of random trigonometric series
in $n$~variables (specifically, Theorem~4 of Chapter~6
of~\cite{Kahane}), there is a constant~$C$ such that for every
collection of complex numbers~$c_\alpha$ and every integer~$M$
greater than~$1$, there is a choice of plus and minus signs for
which the supremum of the modulus of $\sum_{|\alpha|=M}
{\pm}c_\alpha z^\alpha$ in the unit polydisc is no more than $C(n
\sum_{|\alpha|=M}|c_\alpha|^2 \log M)^{1/2}$. We emphasize that
$C$~is independent of the dimension~$n$ and the degree~$M$.

If $r$~is any number less than the Bohr radius~$K_n$, then it
follows from this estimate that
\begin{equation}
\label{eqn:random}
r^M\sum_{|\alpha|=M}|c_\alpha|\le C(n
\sum_{|\alpha|=M}|c_\alpha|^2 \log M)^{1/2}.
\end{equation}
In~(\ref{eqn:random}), we take $c_\alpha$ equal to
$M!/\alpha!$, observing that $\sum_{|\alpha|=M}(M!/\alpha!)=n^M$ on
the left-hand side. On the right-hand side, we crudely estimate
$\sum_{|\alpha|=M} (M!/\alpha!)^2\le
M!\,\sum_{|\alpha|=M}(M!/\alpha!) = M!\,n^M$. Consequently, we
obtain
\begin{equation}
\label{eqn:Mth}
r^M\le C n^{(1-M)/2} (M!\,\log M)^{1/2}.
\end{equation}
Taking $M$th roots
in~(\ref{eqn:Mth}) and choosing~$M$ to be an integer close to
$\log n$, we obtain $K_n< C\sqrt{\log n}/\sqrt{n}$ (with a new
constant~$C$). 

We now make a rough estimate for~$C$. Since $K_n$ never exceeds
$1/3$, the upper bound for~$K_n$ stated in the theorem is
interesting only when $n\ge 189$. If we take $M$~to be the next
integer above $\log n$, then $M\ge 6$ for such~$n$. Now Theorem~1
of Chapter~6 of~\cite{Kahane} with $\kappa=3$ and $\rho=(2\pi
M^2)^n$ shows the existence of a random homogeneous polynomial
$\sum_{|\alpha| =M}{\pm} (M!/\alpha!)z^\alpha$ whose modulus is
bounded on the unit polydisc in dimension~$n$ by
$3\bigl(\sum_{|\alpha|=M}(M!/\alpha!)^2 n\log(6^{1/n}2\pi
M^2)\bigr)^{1/2}$. For the values of~$M$ of interest,
$6^{1/n}2\pi<M^2$, so since $\log M<M$, this upper bound is less
than $6(M!\,n^{M+1}M)^{1/2}$. Arguing as above, we divide by
$n^M$ and take the $M$th root to estimate the Bohr radius. Using
that $n^{1/M}< n^{1/\log n}=e$, while $M!M\le M^M 5!/6^4$ and
$M<(\log n)(1+(\log{189})^{-1})$, we find $K_n<2\sqrt{\log
n}/\sqrt{n}$ as claimed.
\end{proof}

\begin{remark}
Our method seems unlikely to yield the exact value of the Bohr
radius~$K_n$, but small improvements over Theorem~\ref{bohr} are
feasible. For example, the Schwarz lemma for polydiscs (see
Lemma~7.5.6 of~\cite{Rudin}) implies that
$\sum_{|\alpha|=1}|c_\alpha|\le (1-|c_0|^2)$; and if $z$~is in
the polydisc of radius~$r$, then $\sum_{|\alpha|=k}|z^\alpha|^2 \le
r^{2k} \sum_{|\alpha|=k}1 = r^{2k}\tbinom{n+k-1}{k}$. Hence we
find in place of equation~(\ref{eqn:main}) that the Bohr radius
is no smaller than the solution~$r$ of the equation
$r+\sum_{k=2}^\infty r^k \tbinom{n+k-1}{k}^{1/2}=1/2$.

If we put $r=q/\sqrt{n}$ and optimize the value of~$q$, we find
(for example) that $K_n$~exceeds $(2/5)/\sqrt{n}$
when $n>1$, and $K_n$~exceeds $(1/2)/\sqrt{n}$ for 
sufficiently large~$n$.

Also, it is easy to deduce from inequality~(\ref{eqn:Mth}) the
asymptotic upper bound $\limsup_{n\to\infty} K_n\sqrt{n/\log n}\le1$.
\end{remark}

\begin{remark}
We thank Professor Henry Helson for bringing to our attention a
paper of Bohnenblust and Hille \cite{Hille} that constructs
special $M$-linear forms with unimodular coefficients, the forms
being bounded on the unit polydisc in dimension~$n$ by
$n^{(M+1)/2}$. One might hope similarly to construct
\emph{symmetric} $M$-linear forms with unimodular coefficients,
the forms admitting some weaker upper bound $C^M n^{(M+1)/2}$ for
a constant~$C$.  The homogeneous polynomial associated to such a
symmetric form could then be used in our argument above to
eliminate the logarithmic factor in the upper bound for the Bohr
radius. However, we can prove that no such symmetric $M$-linear
form can exist.

In fact, the homogeneous polynomial associated to a symmetric
$M$-linear form with unimodular coefficients can be written
$\sum_{|\alpha|=M}c_\alpha z^\alpha$, where
$|c_\alpha|=M!/\alpha!$. The supremum of such a polynomial
dominates its $L^2$~norm on the unit torus, namely
$(\sum_{|\alpha|=M}(M!/\alpha!)^2)^{1/2}$, which in turn exceeds
$\tbinom{M+n-1}{M}^{-1/2}n^M$ by the Cauchy-Schwarz inequality
applied to $\sum_{|\alpha|=M}[1\cdot (M!/\alpha!)]$. When $n$~is
fixed, this lower bound grows faster than $n^M/M^{n/2}$ as
$M\to\infty$. Hence the modulus of the polynomial cannot admit an
upper bound on the $n$-dimensional unit polydisc of the form $C^M
n^{(M+1)/2}$ with $C$~independent of $n$ and~$M$.\label{remark}
\end{remark}

\begin{remark}
In~\cite{Dineen one} (see also \cite{Dineen two}), S.~Dineen and
R.~M.~Timoney state a result (Theorem~3.2) that, specialized to
polydiscs whose $n$~radii are all equal, says $K_n\le 2/\sqrt n$,
a better asymptotic upper bound than the one in our
Theorem~\ref{bohr}.  However, the supporting Lemma~3.3
in~\cite{Dineen one} is false: it claims the existence of a
symmetric $M$-linear form with unimodular coefficients and
an upper bound $2^{M+1}[Mn^{M+1}\log(1+4M)]^{1/2}+1$
on the unit polydisc in dimension~$n$; but Remark~\ref{remark}
shows that no such form can exist when the dimension~$n$ is
large.  The error in \cite{Dineen one} results from a mistaken
assumption that the constants in the estimates of A.~M.~Mantero
and A.~Tonge \cite{Mantero Tonge} for norms of random tensors
carry over unchanged to the case of \emph{symmetric}
tensors. Adjustments to the proof in \cite{Mantero Tonge} are
required in the symmetric case because the random components of
the tensors are no longer \emph{independent} random variables:
see~\cite{Varopoulos} for a clear exposition of random symmetric
tensors of order three.

\end{remark}

\end{document}